\documentclass[oneside,a4paper]{amsart}

\usepackage{gauss-bonnet}
\hypersetup{
pdftitle={Self-crossing geodesics},
pdfauthor={Anton Petrunin}
}

\begin{document}

\title{Self-crossing geodesics}
\author{Anton Petrunin}
\maketitle

\section{Introduction}

The rubber band on the picture is pulled around a pebble,
and it crosses itself at several points.
\begin{figure}[!ht]
\hfill
\begin{minipage}{.56\textwidth}
\centering
\includegraphics[width=\textwidth]{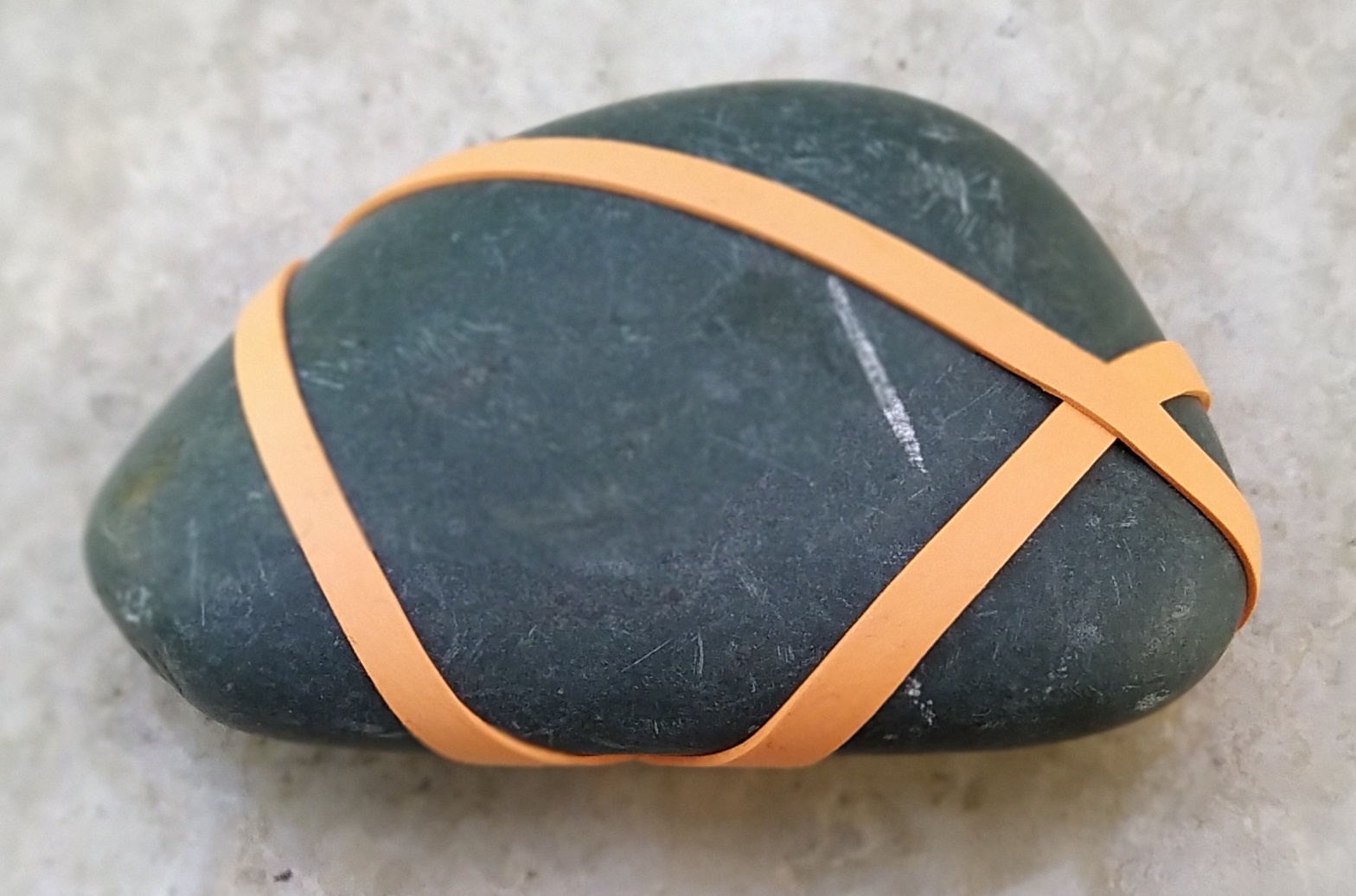}
\end{minipage}
\hfill
\begin{minipage}{.30\textwidth}
\centering
\includegraphics{mppics/pic-50}
\end{minipage}
\hfill
\end{figure}
The combinatorics of self-crossings can be described by a closed plane curve --- it is the rubber band in a parametrization of the surface with one point removed.
For example, if you could turn the pebble around you would see that the self-crossings are described by the plane curve on the right diagram.

We assume that the surface of the pebble is strongly convex, smooth, and frictionless;
in this case, the rubber band models a closed geodesic.
Suppose that we are interested in possible patterns of self-crossings; more precisely:

\medskip

\emph{What are the possible combinatoric types of self-crossings of a closed geodesic on a strongly convex smooth closed surface?}

\begin{figure}[ht!]
\begin{center}
\includegraphics{mppics/pic-55}
\end{center}
\end{figure}

Consider the six possible patterns with three double crossings.
Configurations 1, 2, 3, and 4 can be realized as mirror-symmetric geodesics on mirror-symmetric surfaces; the projections on the plane of symmetry are sketched.
\begin{figure}[ht!]
\begin{center}
\includegraphics{mppics/pic-100}
\end{center}
\end{figure}

Further, we will discuss \emph{forbidden} configurations;
that is, configurations that cannot appear for closed geodesic.
These are configurations 5 and 6.

This question is a good exercise --- it could be explained to anyone, but an answer requires a considerable part of the theory.
The reader is welcome to check that there are no forbidden patterns with less than 3 double crossings and
try the cases with more self-crossings (check \cite[Figures 15--17]{arnold}).
By the way, \emph{is there an algorithm for solving the general case?} 
Our question is closely related to the so-called \emph{flat knot types} of geodesics;
see \cite{angenent} and the references therein.

In what follows, we discuss the Gauss--Bonnet formula as well as the Alexandrov--Toponogov theorem and apply them to forbid configurations 5 and 6.
These theorems are covered in our textbook \cite{petrunin-zamora} which I like, altho they are treated in plenty of other places as well.

\section*{Gauss--Bonnet and no 5}

Suppose that $\Delta$ is an $n$-gon with geodesic sides in a surface $\Sigma$.
Recall that by the \emph{Gauss--Bonnet formula} the sum of the external angles of $\Delta$ equals
\[2\cdot\pi\cdot\chi(\Delta)-\int\limits_\Delta K,\]
where $\chi(\Delta)$ denotes the \emph{Euler characteristic} of $\Delta$ and $K$, the Gauss curvature of~$\Sigma$.

Further, we assume that $\Sigma$ is a closed strongly convex surface. In this case,
\begin{itemize}
 \item $\Sigma$ has strictly positive Gauss curvature;
 \item $\Sigma$ is homeomorphic to the sphere and therefore $\chi(\Sigma)=2$;
 \item 
 $\Delta$ is homeomorphic to the disc and therefore $\chi(\Delta)=1$.
\end{itemize}
It follows that the sum of the internal angles of $\Delta$ is lager than $(n\z-2)\cdot\pi$.
In particular, if $\Delta$ is a triangle with angles $\alpha$, $\beta$, and $\gamma$, then
\[\alpha+\beta+\gamma>\pi.\leqno({*})\]
The Gauss--Bonnet formula can be applied to the whole surface; it implies that
the integral of Gauss curvature along  $\Sigma$ is exactly $4\cdot\pi$.

{

\begin{wrapfigure}{r}{33 mm}
\vskip-8mm
\centering
\includegraphics{mppics/pic-3}
\end{wrapfigure}

\parit{No 5 is forbidden.}
Suppose there is a geodesic with self-crossings as on the diagram;
it divides the surface $\Sigma$ into one triangle, say $\Delta$, one hexagon, and three monogons.
Denote by $\alpha$, $\beta$, and $\gamma$ the internal angles of $\Delta$.

Note that three monogons have internal angles $\alpha$, $\beta$, and~$\gamma$.
By Gauss--Bonnet, the integral of Gauss curvature along each monogon is 
$\pi\z+\alpha$, $\pi\z+\beta$, and $\pi\z+\gamma$ 
respectively.
By $({*})$ the integral of Gauss curvature along the three monogons exceeds $4\cdot \pi$.
But $4\cdot \pi$ is the integral of Gauss curvature along the \emph{whole} surface  --- a contradiction.

}

\section*{Alexandrov--Toponogov and no 6}

Let $\Delta$ be a geodesic triangle with angles $\alpha$, $\beta$, and $\gamma$ on the surface $\Sigma$.
Assume that the sides of $\Delta$ are length-minimizing among the curves \emph{in} $\Delta$ with the same endpoints, then the inequality $({*})$ can be made more exact.

Namely consider the so-called \emph{model triangle} $\tilde\Delta$ of $\Delta$; that is, $\tilde\Delta$ is a plane triangle with equal corresponding sides.
Since the sides are length-minimizing, they satisfy the triangle inequality; therefore the model triangle is defined.

Denote by $\tilde \alpha$, $\tilde \beta$ and $\tilde \gamma$ the angles of $\tilde\Delta$ respectively.
Then 
\[
\alpha> \tilde \alpha,
\qquad
\beta> \tilde \beta,
\qquad
\text{and}
\qquad
\gamma> \tilde \gamma.
\leqno({*}{*})
\]
Since $\tilde\alpha+\tilde\beta+\tilde\gamma=\pi$, this inequality implies $({*})$.

The inequality $({*}{*})$ easily follows from the proof of the \emph{Alexandrov--Toponogov theorem}.
The latter implies that $({*}{*})$ holds for triangles with length-minimizing sides in the \emph{whole} surface.
The proof is left as an 
exercise for those familiar with the Alexandrov--Toponogov theorem; others may simply accept it as true.

\parit{No 6 is forbidden.}
Suppose that such a geodesic $\xi$ exists;
assume that its arcs and angles are labeled as in the leftmost part of the diagram below.
Applying the Gauss--Bonnet formula to the quadrangle and pentagon that $\xi$ cuts from the surface, we get that
\[2\cdot\alpha<\beta+\gamma,
 \quad
2\cdot\beta+2\cdot \gamma<\pi+\alpha,
\quad\text{and therefore} \quad \alpha <\tfrac \pi 3.\leqno(\asterism)\]

\begin{figure}[!ht]
\vskip-1mm
\begin{minipage}{.22\textwidth}
\centering
\includegraphics{mppics/pic-472}
\end{minipage}
\hfill
\begin{minipage}{.35\textwidth}
\centering
\includegraphics{mppics/pic-473}
\end{minipage}
\hfill
\begin{minipage}{.35\textwidth}
\centering
\includegraphics{mppics/pic-474}
\end{minipage}
\vskip-1mm
\end{figure}

Consider the part of $\xi$ without the arc labeled by~$a$.
It cuts from the surface a pentagon $\Delta$ with sides and angles shown in the middle part of the diagram.

Let us add additional vertices on the sides of $\Delta$ so that each side becomes length-minimizing in $\Delta$.
Choose a vertex of $\Delta$ and join it by shortest paths in $\Delta$ to every other vertex.
Consider a model triangle for each triangle in the obtained subdivision of $\Delta$;
the model triangles lie in the plane and we suppose that they share sides as in $\Delta$.
By the comparison inequality $({*}{*})$, the angles of the model triangles do not exceed the corresponding angles of the original triangle.
Therefore, the model triangles form a convex plane polygon, say $\tilde\Delta$,
such that
\begin{itemize}
\item The five angles of $\tilde\Delta$ that correspond to the angles of $\Delta$ do not exceed those.
\item Each side of $\tilde\Delta$ equals to the corresponding small side of $\Delta$.
\end{itemize}
It remains to show that no convex plane polygon meets these two conditions.

\begin{wrapfigure}{r}{47 mm}
\vskip-2mm
\centering
\includegraphics{mppics/pic-475}
\bigskip
\includegraphics{mppics/pic-500}
\vskip1mm
\end{wrapfigure}

Let us orient the sides $\tilde\Delta$ counterclockwise;
denote the obtained vectors by $s_1,\dots,s_k$.
Note that the vectors $s_i$ point in the complement of white sectors shown below with angles marked.
The sum of the magnitudes of the vectors in each black sector is also marked (each black sector corresponds to a side of $\Delta$).

By $(\asterism)$ we can choose a vector $r$ as on the diagram, so that $\phi\z>\tfrac{\pi-\beta}2$ and $\psi>\tfrac{\pi-\gamma}2$.
Note that for any unit vectors $u$, $v$, $v'$, $w$, and $w'$ in the marked black sectors,
we have
\[\begin{aligned}
\langle r,u\rangle<0,\quad \langle r,v\rangle<0,\quad \langle r,w\rangle&<0,
\\
\langle r,v\rangle+\langle r,v'\rangle<0,\quad \langle r,w\rangle+\langle r,w'\rangle&<0.
 \end{aligned}
\qquad\qquad\qquad\qquad\qquad\qquad\qquad\]
It follows that 
\[
\langle r,s_1\rangle+\dots+\langle r,s_k\rangle<0.\qquad\qquad\qquad\qquad\qquad\qquad
\]
But the vectors $s_i$ circumambulate $\tilde\Delta$;
so, the sum has to vanish
--- a contradiction.

{\small \parbf{Acknowledgments.}
I want to thank Arseniy Akopyan, Maxim Arnold, Serge Tabachnikov, David Kramer, and Sergio Zamora Barrera for their help.
This work was partially supported by the NSF grant DMS-2005279, and Simons Foundation grant no 584781.}

{\sloppy
\printbibliography
\fussy
}

\end{document}